\documentclass[oneside, 12pt]{amsart}
\usepackage[utf8]{inputenc}

\usepackage{amsfonts}
\usepackage{amsmath}
\usepackage{amssymb}
\usepackage{amsthm}
\usepackage[a4paper]{geometry}
\usepackage{mathrsfs} 
\usepackage[dvipsnames]{xcolor}

\usepackage[pagebackref=true,colorlinks=true, linkcolor=RedViolet, citecolor=RedViolet]{hyperref}
\renewcommand*\backref[1]{\ifx#1\relax \else (Cited on #1) \fi}

\usepackage{appendix}
\usepackage{enumerate}

\theoremstyle{plain}
\newtheorem{theorem}{Theorem}
\newtheorem{proposition}[theorem]{Proposition}
\newtheorem{lemma}[theorem]{Lemma}
\newtheorem{corollary}[theorem]{Corollary}
\newtheorem{definition}[theorem]{Definition}
\newtheorem{remark}[theorem]{Remark}

\theoremstyle{definition}

\numberwithin{theorem}{section}
\numberwithin{equation}{section} 

\newcommand{\set}[1]{\left\{ #1 \right\}}
\newcommand{\norm}[1]{\left \lVert  #1 \right \rVert}
\newcommand{\abs}[1]{\left\lvert #1 \right\rvert}

\newcommand{\vertiii}[1]{{\left\vert\kern-0.25ex\left\vert\kern-0.25ex\left\vert #1 
    \right\vert\kern-0.25ex\right\vert\kern-0.25ex\right\vert}}

\newcommand{\Z}{\mathbb{Z}}

\newcommand{\R}{\mathbb{R}}
\newcommand{\N}{\mathbb{N}}

\newcommand{\calF}{\mathcal{F}}

\newcommand{\calM}{\mathcal{M}}

\title[Long-time behaviour of IPS]{On the long-time behaviour of reversible interacting particle systems in one and two dimensions}
\author{Benedikt Jahnel}
\address{Institut f\"ur Mathematische Stochastik, Technische Universit\"at Braunschweig, Universit\"atsplatz 2,
38106 Braunschweig, Germany \& Weierstrass Institute for Applied Analysis and Stochastics\\
Mohrenstraße 39\\
10117 Berlin\\
Germany}
\email{benedikt.jahnel@tu-braunschweig.de}

\author{Jonas Köppl}
\address{Weierstrass Institute for Applied Analysis and Stochastics\\
Mohrenstraße 39\\
10117 Berlin\\
Germany}
\email{jonas.koeppl@wias-berlin.de}

\date{\today}
\keywords{Interacting particle systems, gibbs measures, relative entropy, attractor}
\subjclass{Primary 82C22; Secondary 60K35}

\begin{document}

\begin{abstract}
  \noindent
  By refining Holley's free energy technique, we show that, under quite general assumptions on the dynamics, a (possibly non-translation-invariant) interacting particle system in one or two spatial dimensions cannot exhibit time-periodic behaviour if the dynamics admits a reversible Gibbs measure. 
  This is the first result that makes the physical intuition rigorous that time-periodic behaviour can only happen in driven, i.e., non-reversible systems.
\end{abstract}
\maketitle

\section{Introduction and motivation}
One major part of the literature on interacting particle systems deals with the study of their long-time behaviour, in particular the convergence to time-stationary measures and the question of ergodicity. For continuous-time Markov processes on finite state spaces, this is a rather simple question and very well-understood, but for interacting particle systems on $\Z^d$ this is much more subtle. For example, in \cite{jahnel_class_2014} it was shown that for dimension $d\geq 3$, there are non-degenerate irreducible systems with a unique time-stationary measure that fail to be ergodic due to the existence of a periodic orbit in the associated measure-valued dynamics. This is a type of complex behaviour that simply does not occur for irreducible continuous-time Markov processes on finite state spaces.

While a classical result by Mountford, see \cite{mountford_coupling_1995}, shows that no such system can exist in one spatial dimension, it is an open problem whether a two-dimensional interacting particle system can exhibit non-trivial time-periodic behaviour. By now, there is a variety of mean-field systems that have been shown to exhibit time-periodic behaviour, see \cite{collet_rhythmic_2016} or \cite{dai_pra_curie-weiss_2013} for recent results in this spirit. Probably the most famous and classical example for such behaviour is the Kuramoto model which has been studied very successfully, see for example \cite{acebron_kuramoto_2005} and \cite{giacomin_global_2012}. However, according to numerical experiments all of the known examples seem not to exhibit periodic behaviour in two dimensions. 

In this article, we want to investigate the possible long-time behaviour of interacting particle systems in one and two spatial dimensions under the assumption that they admit at least one reversible measure. We show that the presence of a single such fixed point narrows down the possible types of behaviours dramatically and at least morally describes the long-term behaviour completely. 

In particular, we show that the presence of at least one reversible fixed point of the measure-valued dynamics implies that there can be no periodic orbits. This provides a first partial answer to an open question which asks for the existence of interacting particle systems that exhibit time-periodic behaviour in two spatial dimensions, see e.g.~\cite[Section 1.2]{maes_rotating_2011} or \cite[Section 1.8]{swart_course_2022}.

It is often argued that the presence of effects that break the time-symmetry is needed in order for an interacting particle system to exhibit time-periodic behaviour. To our knowledge, our proof gives the first rigorous mathematical justification for using this heuristic for spatially extended systems in continuous time. 

We expect that a similar result should also hold in dimensions $d \geq 3$, but the method of proof breaks down and one needs to proceed differently. It is also not clear if one can or cannot extend this method of proof to non-reversible systems. 

For our proofs, we use a Lyapunov-functional approach which was pioneered in the context of lattice systems by Holley in \cite{holley_free_1971} and later extended to more general and even non-reversible systems in \cite{higuchi_results_1975,kunsch_non_1984,maes_Gibbsian_1994,jahnel_attractor_2019,jahnel_dynamical_2022}.
Roughly speaking, these results show that, if an interacting particle system admits a shift-invariant time-stationary Gibbs measure, then all other shift-invariant time-stationary measures for the dynamics are also necessarily Gibbs with respect to the same specification.
However, all of these papers heavily rely on certain subadditivity properties and therefore only apply to \textit{translation-invariant} measures and cannot be used to say anything meaningful if one starts in a non-translation-invariant measure. 
The only results on non-translation-invariant measures that were obtained via the free-energy method are contained in \cite{holley_one_1977}. There, the authors were able to show that every time-stationary measure of a stochastic Ising model in one and two spatial dimensions is even reversible, hence a Gibbs measure. 

We extend their results to show that not just every time-stationary measure is reversible, but actually every time-stationary orbit is trivial in the sense that it consists of a single reversible measure. 
Therefore, it is not possible to both admit a reversible Gibbs measure and exhibit time-periodic behaviour.

The main conceptual contribution of this article can therefore be summarised as the realisation that the finite-volume free energy technique from \cite{holley_one_1977} can not only be used for analysing stationary measures but can also be applied to \textit{stationary orbits} and thereby yields much more information about the long-term behaviour of interacting particle systems than previously thought. In particular, we show that free energy techniques can be used to rule out time-periodic behaviour in reversible one and two-dimensional interacting particle systems. 

Note that, since we are working with techniques from \cite{holley_one_1977}, we do \textit{not} need to assume any translation-invariance. This is due to being able to only work with finite-volume relative entropy losses and not taking any density limits, which is only possible in dimensions one and two. 

The concept of relative entropy decay has become a very powerful tool for studying systems of interacting particles, and it connects probability, analysis, and geometry in an intricate way. One particularly fruitful application of relative entropy techniques is in the context of Log-Sobolev inequalities for Markov processes. These inequalities can be used to obtain bounds on the (exponential) speed of convergence to equilibrium. However, these methods are limited to the situation where the time-stationary measure is unique, whereas our method goes beyond this case and is also applicable in the non-uniqueness regime. A pedagogical introduction to Log-Sobolev inequalities in the easier setting of Markov chains on finite state spaces can be found in \cite{diaconis_logarithmic_1996}, while a very general approach can be found in \cite[Chapter 5]{bakry_analysis_2014}. 

Another sub-area of probability where relative-entropy methods have successfully been applied to obtain limit theorems is the derivation of hydrodynamic equations as scaling limits of microscopic models of systems of interacting particles. In this context, the method is used to study the infinite particle limit, with additional rescaling of space and time, and not for long-time asymptotics. An introduction to this method can for example be found in the monograph \cite{kipnis_scaling_1999}. 

\subsection{Structure of the manuscript}
In Section~\ref{section:setting-and-notation} we introduce the framework in which we are working and setup the required notation before we state our main results in Section~\ref{section:main-results}. We then proceed by explaining the structure and key ideas of the main proof in Section~\ref{section:proof-strategy}. After this, we finally start with the main work and carry out the proof of our main results in Sections~\ref{section:holley-stroock-principle} and \ref{section:positive-mass}. In the end, we comment on possible extensions to more general interacting particle systems in Section~\ref{section:generalisations}.

\section{Setting and notation}\label{section:setting-and-notation}
Let $q \in \N$ and consider the configuration space $\Omega := \Omega_0^{\Z^d} = \{1,\dots,q\}^{\Z^d}$, which we will equip with the usual product topology and the corresponding Borel sigma-algebra $\calF$. For $\Lambda \subset \Z^d$ let $\calF_\Lambda$ be the sub-sigma-algebra of $\calF$ that is generated by the open sets in $\Omega_\Lambda := \{1,\dots, q\}^{\Lambda}$. We will use the shorthand notation $\Lambda \Subset \Z^d$ to signify that $\Lambda$ is a finite subset of $\Z^d$. In the following we will often denote for a given configuration $\omega \in \Omega$ by $\omega_\Lambda$ its projection to the volume $\Lambda \subset \Z^d$ and write $\omega_\Lambda \omega_\Delta$ for the configuration on in $\Lambda \cup \Delta$ composed of $\omega_\Lambda$ and $\omega_\Delta$ for disjoint $\Lambda, \Delta \subset \Z^d$. For the special case $\Lambda = \{x\}$ we will also write $x^c = \Z^d \setminus \{x\}$ and $\omega_x\omega_{x^c}$. 
The set of probability measures on $\Omega$ will be denoted by $\calM_1(\Omega)$ and the space of continuous functions by $C(\Omega)$. For a configuration $\eta \in \Omega$ we will denote by $\eta^{x,i}$ the configuration that is equal to $\eta$ everywhere except at the site $x$ where it is equal to $i$. Moreover, for $\Lambda \subset \Z^d$ we will denote the corresponding cylinder sets by 
$
    [\eta_\Lambda] = \{\omega : \ \omega_\Lambda \equiv \eta_\Lambda \}. 
$
Whenever we are taking the probability of such a cylinder event with respect to some measure $\nu \in \calM_1(\Omega)$, we will omit the square brackets and simply write $\nu(\eta_\Lambda)$. 

\subsection{Gibbs measures and interacting particle systems}
\subsubsection{Gibbs measures}
We begin by recalling the definition of a specification. 
\begin{definition}
A \emph{specification} $\gamma = (\gamma_\Lambda)_{\Lambda \Subset \Z^d}$ is a family of probability kernels $\gamma_{\Lambda}$ from $\Omega _{\Lambda^c}$ to $\calM_1(\Omega)$ that additionally satisfies the following properties. 
\begin{enumerate}[i.]
    \item Each $\gamma_{\Lambda}$ is \emph{proper}, i.e., if $\Delta \subset \Lambda^c$, then 
    \begin{align*}
        \gamma_\Lambda(\eta_{\Lambda}\eta_\Delta | \eta_{\Lambda^c}) = \gamma_{\Lambda}(\eta_{\Lambda}|\eta_{\Lambda^c}) \mathbf{1}_{\eta_\Delta}(\eta_{\Lambda^c}). 
    \end{align*}
    \item The probability kernels are \emph{consistent} in the sense that if $\Delta \subset \Lambda \Subset \Z^d$, then 
    \begin{align*}
        \gamma_{\Lambda}(\gamma_{\Delta}(\eta_{\Delta}|\cdot)|\eta_{\Lambda^c}) 
        = 
        \gamma_\Lambda(\eta_\Delta | \eta_{\Lambda^c}). 
    \end{align*}
\end{enumerate}
\end{definition}

An infinite-volume probability measure $\mu$ on $\Omega$ is called a \textit{Gibbs measure} for a specification $\gamma$ if $\mu$ satisfies the so-called \textit{DLR equations}, namely for all $\Lambda \Subset \Z^d$  and $\eta_\Lambda$ we have
\begin{align}\label{dlr-equations}
    \mu(\gamma_\Lambda(\eta_\Lambda | \cdot)) = \mu(\eta_\Lambda). 
\end{align}
We will denote the set of all Gibbs measures for a specification $\gamma$ by $\mathscr{G}(\gamma)$. For the existence and further properties of Gibbs measures one needs to impose some conditions on the specification $\gamma$. One sufficient condition for the existence of a Gibbs measure for a specification $\gamma$ is \textit{quasilocality}, which should be thought of as a continuous dependence on the boundary condition.

\begin{definition}
A specification $\gamma$ is called 
\begin{enumerate}[i.]
    \item \emph{non-null}, if for some $\delta >0$
    \begin{align*}
        \inf_{\eta \in \Omega, x \in \Z^d}\gamma_{x}(\eta_0 | \eta_{x^c}) \geq \delta.
    \end{align*}
    \item \emph{quasilocal}, if for all $\Lambda \Subset \Z^d$
    \begin{align*}
        \lim_{\Delta \uparrow \Z^d}\sup_{\eta, \xi \in \Omega}\abs{\gamma_{\Lambda}(\eta_\Lambda | \eta_{\Delta \setminus \Lambda}\xi_{\Delta^c}) - \gamma_\Lambda(\eta_{\Lambda}|\eta_{\Lambda^c})} = 0. 
    \end{align*}
\end{enumerate}
\end{definition}
We will sometimes consider the probability kernels $\gamma_\Lambda$ as functions $\Omega\to[0,1]$, $ \omega\mapsto\gamma_\Lambda(\omega_\Lambda | \omega_{\Lambda^c})$.
If $\gamma$ is a quasilocal specification, then each such map is uniformly continuous.
For example, specifications defined via a uniformly absolutely summable potential $\Phi=(\Phi_B)_{B \Subset \Z^d}$ are non-null and quasilocal. For more details on Gibbs measures and specifications see~\cite{georgii_gibbs_2011},~\cite[Chapter 6]{friedli_statistical_2017} and~\cite[Chapter 4]{bovier_statistical_2006}. 

\subsubsection{Interacting particle systems}

We will consider time-continuous Markovian dynamics on $\Omega$, namely interacting particle systems characterized by time-homogeneous generators $\mathscr{L}$ with domain $\text{dom}(\mathscr{L})$ and its associated Markov semigroup $(P_t)_{t \geq 0}$. 
For interacting particle systems we adopt the notation and exposition of the classical textbook \cite[Chapter 1]{liggett_interacting_2005}. 
In our setting, the generator $\mathscr{L}$ is given via a collection of  single-site transition rates $c_x(\eta, \xi_x)$, which are continuous in the starting configuration $\eta \in \Omega$. 
These rates can be interpreted as the infinitesimal rate at which the particle at site $x$ switches from the state $\eta_x$ to $\xi_x$, given that the rest of the system is currently in state $\eta_{x^c}$. 
The full dynamics of the interacting particle system is then given as the superposition of these local dynamics, i.e., 
\begin{align*}
    \mathscr{L}f(\eta) = \sum_{x \in \Z^d}\sum_{\xi_x}c_x(\eta, \xi_x)[f(\xi_x \eta_{x^c}) - f(\eta)].
\end{align*}
In \cite[Chapter 1]{liggett_interacting_2005} it is shown that the following two conditions are sufficient to guarantee the well-definedness. 
\begin{enumerate}[\bfseries (L1)]
    \item The rate at which the particle at a particular site changes its spin is uniformly bounded, i.e.,
    \begin{align*}
        \sup_{x \in \Z^d}\sum_{\xi_x}\norm{c_{x}(\cdot, \xi_{x})}_{\infty} < \infty 
    \end{align*}
    \item and the total influence of all other particles on a single particle is uniformly bounded, i.e.,
    \begin{align*}
        \sup_{x \in \Z^d}\sum_{y \neq x}\sum_{\xi_{x}}\delta_y\left(c_{x}(\cdot, \xi_{x})\right) < \infty, 
    \end{align*}
    where 
    \begin{align*}
        \delta_y(f) := \sup_{\eta, \xi\colon \eta_{y^c} = \xi_{y^c}}\abs{f(\eta)-f(\xi)}
    \end{align*}
    is the oscillation of a function $f:\Omega \to \R$ at the site $y \in \Z^d$. 
\end{enumerate}
Under these conditions, one can then show that the operator $\mathscr{L}$, defined as above, is the generator of a well-defined Markov process and that a core of $\mathscr{L}$ is given by the space of functions with finite total oscillation, i.e.
\begin{align*}
    D(\Omega) := \Big\{ f \in C(\Omega)\colon \sum_{x \in \Z^d} \delta_x(f) < \infty\Big\}.
\end{align*}

\noindent Let us emphasise briefly that we will not assume translation-invariance of the rates.

\subsection{Relative entropy loss}
For $\mu, \nu \in \calM_1(\Omega)$ and a finite volume $\Lambda \Subset \Z^d$ define the relative entropy of $\nu$ with respect to $\mu$ in $\Lambda$ via 
\begin{align*}
    h_\Lambda(\nu | \mu) :=
    \begin{cases}
    \sum_{\omega_\Lambda \in \Omega_\Lambda}\nu(\omega_{\Lambda})\log\frac{\nu(\omega_{\Lambda})}{\mu(\omega_\Lambda)}, \quad &\text{if } \nu_\Lambda \ll \mu_\Lambda, 
    \\\
    \infty, &\text{else,}
    \end{cases}
\end{align*}
where we use the convention that $0 \log 0 = 0$. Now recall that $(P_t)_{t \geq 0}$ denotes the Markov semigroup corresponding to the Markov generator $\mathscr{L}$. We write $\nu_t:=\nu P_t$ for the time-evolved measure $\nu\in\calM_1(\Omega)$.
The finite-volume relative entropy loss in $\Lambda \Subset \Z^d$ is defined by 
\begin{align*}
    g_\Lambda^{\mathscr{L}}(\nu|\mu) 
    := \frac{d}{dt}\Big \lvert_{t=0}h_{\Lambda}( \nu_t\lvert \mu). 
\end{align*}
Usually, one then works with the density limits of the relative entropy and the relative entropy loss and shows that the latter can still be used as a Lyapunov function for the dynamics. However, the sub-additivity arguments that are used to show that the relative entropy loss density actually exists as a limit are only available for translation-invariant measures. We do not want to make any such assumptions and therefore we will have to instead work with the family of finite-volume relative entropy losses. Note that calling the finite-volume derivatives \textit{loss} is not entirely correct, since they are not necessarily non-positive. However one can show that the positive contributions are of boundary order and vanish in the density limit, see \cite[Lemma~3.10 and Lemma~3.12]{jahnel_dynamical_2022}. 

\subsection{Time-stationary measures, orbits, and the attractor}
If one is interested in the long-term behaviour of an interacting particle system, a natural object to study is the so-called \textit{attractor} of the measure-valued dynamics which is defined as
\begin{align*}
    \mathscr{A} = \set{\nu \in \calM_1(\Omega)\colon \exists \nu_0 \in \calM_1(\Omega) \text{ and } t_n \uparrow \infty \text{ such that } \lim_{n \to \infty}\nu_{t_n} = \nu}. 
\end{align*}
This is the set of all accumulation points of the measure-valued dynamics induced by $\mathscr{L}$. In the language of dynamical systems this is the $\omega$-limit set. This encodes (most of) the dynamically relevant information about the long-time behaviour of the system. 
In this article, we are particularly interested in two subsets of the attractor, namely the \textit{stationary measures} given by 
\begin{align*}
    \mathscr{S} := \set{\nu \in \calM_1(\Omega): \forall s \geq 0: \nu P_s = \nu},
\end{align*}
and the measures which lie on a \textit{stationary orbit}
\begin{align*}
    \mathscr{O} := \set{\nu \in \calM_1(\Omega): \exists T > 0: \nu P_T = \nu}. 
\end{align*}
The relation between these sets can be summarised as follows
\begin{align*}
    \mathscr{S} \subset \mathscr{O} \subset \mathscr{A}. 
\end{align*}
 In general, the first inclusion is strict as can be seen by considering the non-trivial examples constructed in \cite{jahnel_class_2014} and \cite{jahnel_time_periodic_2024} or the (from a probabilistic point of view) trivial example given in \cite[p.12]{liggett_interacting_2005}. 
Historically, most attention has been paid to investigating the set of time-stationary measures and their properties, but not much was known about the behaviour of interacting particle systems outside of this set.
\section{Main results}\label{section:main-results}
\subsection{Assumptions on the transition rates and the specification}
Before we can state our main results, let us introduce some stronger conditions on the specification $\gamma = (\gamma_\Delta)_{\Delta \Subset \Z^d}$ and the rates $(c_x(\cdot, \xi_x))_{x \in \Z^d, \xi_x \in \Omega_0}$ that will turn out to be crucial for our results. 
\medskip

\noindent
\textbf{Conditions for the specification.} 
\begin{enumerate}[\bfseries (S1)]
\item $\gamma$ is quasilocal.
\item $\gamma$ is non-null with constant $\delta >0$. 
\item $\gamma$ satisfies 
\begin{align*}
    \sum_{y \in \Z^d}\abs{y}\sup_{x \in \Z^d}\delta_{x+y}\left(\gamma_x(\cdot) \right)< \infty. 
\end{align*}
\end{enumerate}

\noindent \textbf{Conditions for the rates.} 
\begin{enumerate}[\bfseries (R1)]
    \item  The rate at which the particle at a particular site changes its spin is uniformly bounded, i.e.,
    \begin{align*}
        \sup_{x \in \Z^d}\sum_{\xi_{x}}\norm{c_{x}(\cdot, \xi_{x})}_{\infty} < \infty. 
    \end{align*}
    \item For every $x \in \Z^d$ and $\xi_x \in \Omega_0$ the function 
    \begin{align*}
        \Omega \ni \eta \mapsto c_x(\eta, \xi_x) \in [0,\infty)
    \end{align*}
    is continuous.
    \item The transition rates are bounded away from zero, i.e.,
    \begin{align*}
        \inf_{x\in \Z^d,\, \eta \in \Omega,\, \xi_x \in \Omega_0 } c_x(\eta, \xi_x) \geq \delta >0. 
    \end{align*}
\item We have 
\begin{align*}
    \sum_{y \in \Z^d} \abs{y}\sup_{x \in \Z^d} \delta_{x+y}c_x(\cdot) < \infty,
\end{align*}
where $c_x(\eta) = \sum_{i \neq \eta_x}c_x(\eta,i)$ is the total rate at which the particle at site $x$ changes its state when the system is in configuration $\eta$. 
\end{enumerate}
The conditions $\mathbf{(S3)}$ and $\mathbf{(R4)}$ essentially ensure that the specification and the transition rates are \textit{short-range} and are satisfied if the dependence of $\gamma_x$ and $c_x$ on the spin at site $y \in \Z^d$ decays faster than $\abs{x-y}^{-2d}$.
Let us emphasise again that we do not assume that either the rates or the specification are translation invariant.

Our main result is the following no-go result that essentially states that the presence of a reversible fixed point for the measure-valued dynamics makes it impossible to also possess periodic orbits. The precise formulation is as follows. 
\begin{theorem}\label{theorem:main-result}
Let $d \in \{1,2\}$ and assume that $\mathscr{L}$ is the generator of an interacting particle system that satisfies assumptions $\mathbf{(R1)-(R4)}$ and that it admits a reversible measure $\mu$ that is a Gibbs measure with respect to a specification $\gamma$ that satisfies $\mathbf{(S1)-(S3)}$. 
Then we have that 
\begin{align*}
    \mathscr{O} = \mathscr{S} = \mathscr{G}(\gamma). 
\end{align*}
\end{theorem}

This limit theorem in particular implies the following no-go result that essentially states that the presence of a reversible fixed point for the measure-valued dynamics makes it impossible to also possess periodic orbits. The precise formulation is as follows. 

\begin{corollary}\label{corollary:no-periodic}
Let $d \in \{1,2\}$ and assume that $\mathscr{L}$ and $\gamma$ satisfy the above assumptions and that $(P_t)_{t \geq 0}$ is the Markov semigroup generated by $\mathscr{L}$ with $\mu \in \mathscr{G}(\gamma)$ as reversible fixed point. Then, the measure-valued dynamics given by 
\begin{align*}
    [0,\infty) \times \calM_1(\Omega) \ni (t, \nu) \mapsto \nu P_t \in \calM_1(\Omega)
\end{align*}
does not contain non-trivial time-periodic orbits, i.e., there is no probability measure $\nu \in \calM_1(\Omega)$ such that $(\nu P_t)_{t\geq 0}$ is non-constant and such that there exists a $T>0$ with $\nu_t = \nu_{t +T}$. 
\end{corollary}

\section{Proof strategy}\label{section:proof-strategy}
The proof of Theorem \ref{theorem:main-result} essentially consists of two main steps that we will now explain briefly before we start with the actual mathematics. Let us already point out that apart from the very last step in the proof of Proposition \ref{proposition:time-averaged-holley-stroock} all the technical results in the forthcoming sections hold in any dimension $d \in  \N$ and can be used for future investigations in arbitrary dimensions. 

\subsection{Finite-volume relative entropy loss and Gibbs measures}

The first technical result is the following \textit{time-averaged} version of the results in \cite{holley_one_1977} which also extends the classical results to general finite local state spaces and specifications. 

\begin{proposition}[Time-averaged Holley--Stroock principle]\label{proposition:time-averaged-holley-stroock}
    Assume that $d \in \{1,2\}$, that $\mathscr{L}$ satisfies $\mathbf{(R1)-(R4)}$ and admits a reversible measure $\mu$ which is a Gibbs measure with respect to a specification $\gamma$ that satisfies $\mathbf{(S1)-(S3)}$. If $\nu \in \calM_1(\Omega)$ is a probability measure that satisfies
    \begin{enumerate}[\bfseries \upshape (M1)]
        \item for all $\eta \in \Omega$, $\Lambda \Subset \Z^d$, and $s \geq 0$ it holds that $\nu P_s(\eta_\Lambda) >0$, and 
        \item there exists $T>0$ such that for all $\Lambda \Subset \Z^d$ we have $\int_0^T g^\mathscr{L}_\Lambda(\nu P_s \lvert \mu) ds = 0$,
    \end{enumerate}
    then we have $\nu = \nu P_s$ for all $s\geq 0$ and $\nu \in \mathscr{G} (\gamma)$.
\end{proposition}

The proof of this can be found in Section \ref{section:holley-stroock-principle} and follows the strategy laid out in \cite[Chapter IV.5]{liggett_interacting_2005} in our more general setting but we additionally need to make a distinction between \textit{pointwise estimates}, i.e., for fixed $s \in [0,T]$, and \textit{averaged estimates}.
Let us point out that morally this characterises Gibbs measures as those points at which the finite-volume relative entropy loss vanishes, so it can be interpreted as a dynamical counterpart to the classical Gibbs variational principle, see e.g.~\cite[Theorem 6.82]{friedli_statistical_2017}. However, in the dynamical setting, at least in dimensions one and two, this also holds without assuming translation invariance, which does not work in the static world. 
After establishing this general principle, our main result Theorem \ref{theorem:main-result} will follow by showing that $\mathbf{(M1)-(M2)}$ are satisfied along time-periodic orbits. Condition $\mathbf{(M2)}$ is directly implied by periodicity and the fundamental theorem of calculus, so one only needs to worry about $\mathbf{(M1)}$.

\subsection{The positive-mass property}
If one considers an irreducible continuous-time Markov chain on a finite state space $X$, then it is easy to show that there exist constants $\rho,\tau>0$ such that for any initial distribution $\nu$ and all $x \in X$ and $t\geq \tau > 0$, we have $\nu_t(x) \geq \rho$. In other words, every state has strictly positive mass for any positive time. In the setting of infinite-volume interacting particle systems that are irreducible in a suitable way, something similar should be true, but here it is not as straightforward to see as in the setting of finite state spaces. Proposition~\ref{proposition:positive-mass} makes this intuition precise. A less general but also stronger result was previously derived in \cite{jahnel_attractor_2019}.

To show this positivity property, we will compare our dynamics to an interacting particle system in which all of the sites inside of $\Lambda$ behave independently and flip with the minimal transition rate. We then use the following Girsanov-type formula to compare this finite-volume perturbation with the original dynamics. 

\begin{lemma}\label{lemma:girsanov-transformation}
   Consider an interacting particle system with generator $L$ such that its transition rates $(c_x(\cdot, \cdot))_{x \in \Z^d}$ satisfy assumptions $\mathbf{(L1)-(L2)}$ and are strictly positive. Let $\hat{L}^\Lambda$ be the generator of another interacting particle system with rates $(\hat{c}_x(\cdot, \cdot))_{x \in \Z^d}$ such that the rates of $L$ and $\hat{L}^\Lambda$ agree for sites outside of the finite volume $\Lambda \Subset \Z^d$. 
   Denote the induced path measures on the space of $\Omega$-valued c\'adl\'ag paths $\sigma[0,\tau]$ up to time $\tau >0$ by $\mathbb{Q}_\omega$ respectively $\hat{\mathbb{Q}}_\omega^\Lambda$, where the initial condition $\sigma(0) = \omega$ is deterministic. Then the following Girsanov-type formula holds 
   \begin{align*}
    \frac{d \mathbb{Q}_\omega}{d\hat{\mathbb{Q}}^\Lambda_\omega} ( \sigma[0,\tau])
    =
    \exp\left(-\int_0^\tau \lambda(\sigma(s)) + \sum_{s \in [0,\tau]: \sigma_{\Lambda}(s_-) \neq \sigma_{\Lambda}(s)}\sum_{i \in \Lambda}\log\left(\frac{c_i(\sigma(s_-), \sigma_i(s))}{\hat{c}_i(\sigma(s_-), \sigma_i(s))}\right)\right), 
\end{align*}
where 
\begin{align*}
    \lambda(\eta) :=\sum_{i \in \Lambda}\left(c_i(\eta) - \hat{c}_i(\eta)\right),
\end{align*}
and $c_i(\eta)$ respectively $\hat{c}_i(\eta)$ are the total rates at which we see a flip at site $i$ when we are currently in configuration $\eta$, i.e.,
\begin{align*}
        c_i(\eta) = \sum_{\xi_i \neq \eta_i}c_i(\eta, \xi_i) \quad \text{respectively} \quad \hat{c}_i(\eta) = \sum_{\xi_i \neq \eta_i}\hat{c}_i(\eta, \xi_i). 
\end{align*}
\end{lemma}

In a similar form this appeared in \cite{pra_large_1993} in the context of large deviations for interacting particle systems. 
With this explicit transformation formula at hand, the problem described above can essentially be reduced to controlling the tails of a Poisson random variable and we obtain the following result. 

\begin{proposition}[Positive-mass property]\label{proposition:positive-mass}
Assume that the rates of an interacting particle system satisfy assumptions $\mathbf{(R1)}-\mathbf{(R4)}$. Then, for all $\tau >0$ and $\Lambda \Subset \Z^d$, there exists a constant $C(\tau, \Lambda) >0$ such that, for any starting measure $\nu$ and any time $t \in [\tau, \infty)$, we have
\begin{align*}
    \forall \eta \in \Omega: \quad \nu_t(\eta_\Lambda) \geq C(\tau, \Lambda). 
\end{align*}
In particular, for all subsequential limits $\nu^* = \lim_{n \to \infty}\nu_{t_n}$ with $t_n \uparrow \infty$, we have 
\begin{align*}
    \forall \eta \in \Omega \ \forall \Lambda \Subset \Z^d: \quad \nu^*(\eta_\Lambda) \geq  C(\tau, \Lambda) > 0. 
\end{align*}
\end{proposition}

Note that this implies that condition $\mathbf{(M1)}$ holds along time-periodic orbits.
On an intuitive level, the above result should be interpreted as a somewhat quantitative version of the diffusive nature of the dynamics. Even if we start our process with a point mass $\delta_\omega$ in $\omega \in \Omega$ as initial condition, the distribution of the process at time $t>0$ will already put positive mass on any cylinder set $[\eta_\Lambda]$. The proof of this can be found in Section~\ref{section:positive-mass}.

\section{Proof of the relative entropy loss principle}\label{section:holley-stroock-principle}
\subsection{Characterising reversible measures}
Before we start with characterising reversible measures we state a technical tool that is reminiscent of Lebesgue's differentiation theorem. 
\begin{lemma}[Differentiation lemma]\label{differentiation-lemma}
    Let $\mu$ be a probability measure on $\Omega$ such that we have $\mu(\eta_{\Lambda})>0$ for all $\Lambda \Subset \Z^d$ and $\eta \in \Omega$. 
    Then, for any continuous functions $f: \Omega \to \R$, we have that for all $\eta \in \Omega$
    \begin{align*}
         \lim_{\Lambda \uparrow \Z^d} \frac{1}{\mu(\eta_{\Lambda})}\int_{[\eta_{\Lambda}]}f(\xi)\mu(d\xi) = f(\eta). 
    \end{align*}
    Moreover, if $f$ is uniformly continuous, then the claimed convergence is also uniform in $\eta \in \Omega$. 
\end{lemma}

\begin{proof}
    First note that, for fixed $\Lambda \Subset \Z^d$, the compactness of $\Omega$ implies the trivial inequalities 
    \begin{align}\label{sandwich}
        - \infty < \inf_{\xi: \xi_{\Lambda}=\eta_{\Lambda}}f(\xi) \leq f(\eta) \leq \sup_{\xi: \xi_{\Lambda}=\eta_{\Lambda}}f(\xi) < \infty. 
    \end{align}
    The continuity of $f$ implies that 
    \begin{align*}
        \lim_{\Lambda \uparrow \Z^d}\inf_{\xi: \xi_{\Lambda}=\eta_{\Lambda}}f(\xi) = f(\eta), \quad \lim_{\Lambda \uparrow \Z^d}\sup_{\xi: \xi_{\Lambda}=\eta_{\Lambda}}f(\xi) = f(\eta).
    \end{align*} 
    Combining this with~\eqref{sandwich} and the squeeze theorem (for nets) from real analysis yields
    \begin{align*}
        \lim_{\Lambda \uparrow \Z^d} \frac{1}{\mu(\eta_{\Lambda})}\int_{[\eta_{\Lambda}]}f(\xi)\mu(d\xi) = f(\eta). 
    \end{align*}
    This concludes the proof. 
\end{proof}

With this technical helper at hand, we begin with the following standard result that provides us with alternative formulations of reversibility. These will turn out to be more convenient to work with.  

\begin{proposition}\label{proposition:characterizations-reversibility}
Consider an interacting particle system with generator $\mathscr{L}$ whose rates satisfy assumptions $\mathbf{(L1)}-\mathbf{(L2)}$ and $\mathbf{(R3)}$. Then, for a probability measure $\nu \in \calM_1(\Omega)$, the following conditions are equivalent. 
\begin{enumerate}[i.]
    \item $\nu$ is reversible. 
    \item For all $\Lambda \Subset \Z^d$, $x \in \Z^d$, $i \in \Omega_0$ and $\eta \in \Omega$ it holds that $\nu(\eta_\Lambda) >0$ and
    \begin{align*}
        \int_{[\eta_\Lambda]}c_x(\omega,i)\nu(d\omega)
        =
        \int_{[\eta_\Lambda^{x,i}]}c_x(\omega,\eta_x)\nu(d\omega). 
    \end{align*}
    \item For all $\Lambda \Subset \Z^d$ and $\eta \in \Omega$ we have $\nu(\eta_\Lambda) >0$ and the conditional marginals of $\nu$ satisfy the detailed balance condition, i.e., $\nu$-almost surely
    \begin{align*}
        \forall x \in \Z^d \ \forall \xi_x \in \Omega_0: \nu(\eta_x \lvert \eta_{x^c})c_x(\eta,\xi_x) = \nu(\xi_x \lvert \eta_{x^c})c_x(\xi_x \eta_{x^c},\eta_x). 
    \end{align*}
\end{enumerate}
\end{proposition}

\begin{proof}
    \textit{Ad $i. \Rightarrow ii.$: } As a first step, note that Proposition~\ref{proposition:positive-mass} applies to $\nu$ which yields $\nu(\eta_\Lambda) > 0$. By reversibility of $\nu$ we know that for all $f,g \in D(\Omega)$
    \begin{align*}
        \int_\Omega f(\omega) \mathscr{L}g(\omega) \nu(d\omega) = \int_\Omega g(\omega) \mathscr{L}f(\omega)\nu(d\omega). 
    \end{align*}
    For fixed $\eta \in \Omega$, $\Lambda \Subset \Z^d$, $x \in \Z^d$ and $i=1,\dots,q$ we can apply this to the functions 
    \begin{align*}
        f = \mathbf{1}_{[\eta_\Lambda]}, \quad g=\mathbf{1}_{[\eta_\Lambda^{x,i}]}. 
    \end{align*}
    Then we have
    \begin{align*}
        \int_\Omega f(\omega) \mathscr{L}g(\omega)
        &=
        \sum_{y \in \Z^d}\sum_{j=1}^q \int_\Omega c_y(\omega, j)\left[\mathbf{1}_{[\eta_\Lambda]}(\omega)\mathbf{1}_{[\eta_\Lambda^{x,i}]}(\omega^{y,j}) - \mathbf{1}_{[\eta_\Lambda]}(\omega)\mathbf{1}_{[\eta_\Lambda^{x,i}]}(\omega)\right]\nu(d\omega)
        \\\
        &=
        \int_\Omega c_x(\omega,i)\mathbf{1}_{[\eta_\Lambda]}(\omega)\nu(d\omega). 
    \end{align*}
    On the other hand 
    \begin{align*}
        \int_\Omega g(\omega) \mathscr{L}f(\omega)
        &=
        \sum_{y \in \Z^d}\sum_{j=1}^q \int_\Omega c_y(\omega, j)\left[\mathbf{1}_{[\eta_\Lambda]}(\omega^{y,j})\mathbf{1}_{[\eta_\Lambda^{x,i}]}(\omega) - \mathbf{1}_{[\eta_\Lambda]}(\omega)\mathbf{1}_{[\eta_\Lambda^{x,i}]}(\omega)\right]\nu(d\omega)
        \\\
        &=
        \int_\Omega c_x(\omega,\eta_x)\mathbf{1}_{[\eta_\Lambda^{x,i}]}(\omega)\nu(d\omega).
    \end{align*}
    So the assumed reversibility of $\nu$ implies
    \begin{align*}
         \int_\Omega c_x(\omega,i)\mathbf{1}_{[\eta_\Lambda]}(\omega)\nu(d\omega)
        =
        \int_\Omega c_x(\omega,\eta_x)\mathbf{1}_{[\eta_\Lambda^{x,i}]}(\omega)\nu(d\omega).
    \end{align*}
    \medskip 

    \noindent 
    \textit{Ad $ii. \Rightarrow iii.$: } We have 
    \begin{align*}
    \int_{[\eta_\Lambda]}c_x(\omega,i)\nu(d\omega) = \int_{[\eta_\Lambda^{x,i}]}c_x(\omega, \eta_x)\nu(d\omega). 
    \end{align*}
    So in particular we get 
    \begin{align*}
    \frac{\nu(\eta_\Lambda)}{\nu(\eta_\Lambda)}\int_{[\eta_\Lambda]}c_x(\omega,i)\nu(d\omega)
    =
    \frac{\nu(\eta_\Lambda^{x,i})}{\nu(\eta_\Lambda^{x,i})}\int_{[\eta_\Lambda^{x,i}]}c_x(\omega, \eta_x)\nu(d\omega). 
    \end{align*}
    Now let us rearrange this to get 
    \begin{align*}
    \frac{\nu(\eta_\Lambda)}{\nu(\eta_{\Lambda}^{x,i})}
    =
    \frac{\nu(\eta_\Lambda^{x,i})^{-1}\int_{[\eta_\Lambda^{x,i}]}c_x(\omega,\eta_x)\nu(d\omega)}{\nu(\eta_\Lambda)^{-1}\int_{[\eta_\Lambda]c_x(\omega,i)\nu(d\omega)}}. 
    \end{align*}
    By Lemma \ref{differentiation-lemma} the right hand side converges to 
    \begin{align*}
    \frac{c_x(\eta^{x,i}, \eta_x)}{c_x(\eta,i)}
    \end{align*}
    and by martingale convergence we see that for $\nu$-almost every $\eta \in \Omega$ as $\Lambda \uparrow \Z^d$
    \begin{align*}
    \frac{\nu(\eta_\Lambda)}{\nu(\eta_{\Lambda}^{x,i})}
    =
    \frac{\nu(\eta_x \lvert \eta_{\Lambda \setminus x})}{\nu(i \lvert\eta_{\Lambda \setminus x})}
    \to 
    \frac{\nu(\eta_x \lvert \eta_{x^c})}{\nu(i \lvert \eta_{x^c})}.
    \end{align*}
    After rearranging this is simply the detailed balance equation. 
    \medskip 

    \noindent 
    \textit{Ad $iii. \Rightarrow i.$: } Here it suffices to show that, for all local functions $f,g : \Omega \to \R$, it holds that 
    \begin{align*}
        \int_\Omega f(\omega) \mathscr{L}g(\omega) \nu(d\omega) = \int_\Omega g(\omega) \mathscr{L}f(\omega) \nu(d\omega). 
    \end{align*}
    To do this, one can proceed exactly as in the proof of \cite[Lemma 3.1]{jahnel_dynamical_2022}. 
\end{proof}

Let us briefly summarise the main implications of Proposition~\ref{proposition:characterizations-reversibility} for what is to come. On the one hand, since we assume that our interacting particle system admits a reversible measure $\mu$ that is a Gibbs measure with respect to a non-null specification $\gamma$, we see that $\mu$-almost surely 
\begin{align*}
    c_x(\eta, i) \gamma_x(\eta_x \lvert \eta_{x^c}) = c_x(\eta^{x,i}, \eta_x) \gamma_x(i \lvert \eta_{x^c}). 
\end{align*}
In particular, every other measure $\mu' \in \mathscr{G}(\gamma)$ is also reversible for our process.
On the other hand, the above characterisation tells us that, in order to show that a measure $\nu$ is reversible, we actually only need to show that it satisfies $\nu(\eta_\Lambda) > 0$ for all $\Lambda \Subset \Z^d$ and $\eta \in \Omega$ and that, for all $x \in \Z^d$ and $i \in \Omega_0$, it holds that 
\begin{align*}
    \int_{[\eta_\Lambda]}c_x(\omega,i)\nu(d\omega)
        =
    \int_{[\eta_\Lambda^{x,i}]}c_x(\omega,\eta_x)\nu(d\omega). 
\end{align*}
So our goal will be to show that if $\nu \in \calM_1(\Omega)$ is such that $g^\mathscr{L}_\Lambda(\nu \lvert \mu) = 0$ for all $\Lambda$, then it necessarily satisfies this equation. 
\subsection{Finite-volume relative entropy loss}
We now proceed by obtaining a more convenient representation of the relative entropy loss in the finite cube $\Lambda_n = [-n,n]^d$. Recall that the relative entropy loss in the finite volume $\Lambda_n$ is defined by 
\begin{align*}
    g^n_{\mathscr{L}}(\nu|\mu) 
    := \frac{d}{dt}\Big \lvert_{t=0}h_{\Lambda_n}( \nu_t\lvert \mu), \quad \nu \in \calM_1(\Omega). 
\end{align*}

The first step in our proof of Proposition ~\ref{proposition:time-averaged-holley-stroock} is the following very convenient rewriting of the relative entropy loss in a finite volume $\Lambda$. This first appeared in \cite{moulin_ollagnier_free_1977} and is really the backbone of the whole proof technique. 
We will come back to the importance of this representation in Section~\ref{section:generalisations}, where we discuss possible directions for future generalisations. 

\begin{lemma}\label{lemma:finite-volume-relative-entropy-loss}
    For $n\in \N$ and $\nu \in \calM_1(\Omega)$ we have 
    \begin{align*}
        2 g^n_{\mathscr{L}}(\nu|\mu) 
        =
    -&\sum_{\eta_{\Lambda_n}}\sum_{x \in \Lambda_n}\sum_{i\neq \eta_x}
    \left[\Gamma_n^\nu(x,\eta_x,\eta^{x,i})-\Gamma_n^\nu(x,i,\eta) \right]
    \log\left(\frac{\Gamma_n^\nu(x, \eta_x, \eta^{x,i})}{\Gamma_n^\nu(x,i,\eta)}\right)
    \\\
    +
    &\sum_{\eta_{\Lambda_n}}\sum_{x \in \Lambda_n}\sum_{i\neq \eta_x}
    \left[\Gamma_n^\nu(x,\eta_x,\eta^{x,i})-\Gamma_n^\nu(x,i,\eta) \right]
    \\\
    &\quad 
    \times \left\{
        \log\left(\frac{\nu(\eta_{\Lambda_n})}{{\Gamma_n^\nu(x,i,\eta)}}\right)
        -
        \log\left(\frac{\nu(\eta^{x,i}_{\Lambda_n})}{\Gamma_n^\nu(x,\eta_x, \eta^{x,i})}\right)
        -
        \log\left(\frac{\mu(\eta_{\Lambda_n})}{\mu(\eta_{\Lambda_n}^{x,i})}\right)
        \right\}, 
    \end{align*}
    where we introduce the notation 
    \begin{align*}
    \Gamma^\nu_n(x,j,\eta) := \int_{[\eta_{\Lambda_n}]}c_x(\omega, j)\nu(d\omega). 
    \end{align*}
\end{lemma}
Before we start with the proof, note that by Proposition \ref{proposition:characterizations-reversibility}, our goal will be to show that, for all $x \in \Lambda_n$ and $i\in \Omega_0$, we have $\Gamma_n^\nu(x,i,\eta) = \Gamma_n^\nu(x,\eta_x, \eta^{x,i})$. 
\begin{proof}
    This can be seen by a direct calculation using the definition of the generator and some subsequent algebraic manipulations. We have 
    \begin{align*}
        g^n_{\mathscr{L}}(\nu\lvert \mu) 
        &=
        \sum_{\eta_{\Lambda_n}}
            \nu(\mathscr{L}\mathbf{1}_{\eta_{\Lambda_n}})
            \log
            \left(
                \frac{\nu(\eta_{\Lambda_n})}{\mu(\Lambda_n)}
            \right)
        \\\
        &=
        \sum_{\eta_{\Lambda_n}}\sum_{x \in \Lambda_n}\sum_{j=1}^q
            \int \nu(d\omega)
                c_x(\omega, j)
                \left[
                    \mathbf{1}_{\eta_{\Lambda_n}}(\omega^{x,j}) 
                    -
                    \mathbf{1}_{\eta_{\Lambda_n}}(\omega)
                    \right]
            \log
            \left(
                \frac{\nu(\eta_{\Lambda_n})}{\mu(\eta_{\Lambda_n})}
            \right).
    \end{align*}
    Now let us rewrite this a little bit to bring it into the nicer form
    \begin{align*}
        g^n_{\mathscr{L}}(\nu|\mu) 
        &=
        \sum_{\eta_{\Lambda_n}}\sum_{x \in \Lambda_n}\sum_{i\neq \eta_x}
        \left[\int_{[\eta_{\Lambda_n}^{x,i}]}c_x(\omega, \eta_x)\nu(d\omega) - \int_{[\eta_{\Lambda_n}]}c_x(\omega,i)\nu(d\omega)\right]\log\left(\frac{\nu(\eta_{\Lambda_n})}{\mu(\eta_{\Lambda_n})}\right).
    \end{align*}
     With the notation introduced above this can be written as
    \begin{align*}
        g^n_{\mathscr{L}}(\nu|\mu) 
        &=
        \sum_{\eta_{\Lambda_n}}\sum_{x \in \Lambda_n}\sum_{i\neq \eta_x}
        \left[\Gamma_n^\nu(x,\eta_x,\eta^{x,i})-\Gamma_n^\nu(x,i,\eta) \right]
        \log\left(\frac{\nu(\eta_{\Lambda_n})}{\mu(\eta_{\Lambda_n})}\right)
        \\\
        &= 
        \frac{1}{2} \sum_{\eta_{\Lambda_n}}\sum_{x \in \Lambda_n}\sum_{i\neq \eta_x}
        \left[\Gamma_n^\nu(x,\eta_x,\eta^{x,i})-\Gamma_n^\nu(x,i,\eta) \right]
        \log\left(\frac{\nu(\eta_{\Lambda_n})}{\mu(\eta_{\Lambda_n})}\frac{\mu(\eta_{\Lambda_n}^{x,i})}{\nu(\eta_{\Lambda_n}^{x,i})}\right)
        \\\
        &= 
        \frac{1}{2} \sum_{\eta_{\Lambda_n}}\sum_{x \in \Lambda_n}\sum_{i\neq \eta_x}
        \left[\Gamma_n^\nu(x,\eta_x,\eta^{x,i})-\Gamma_n^\nu(x,i,\eta) \right]
        \log\left(\frac{\nu(\eta_{\Lambda_n})}{\nu(\eta_{\Lambda_n}^{x,i})}\right)
        \\\
        &\ -
        \frac{1}{2}\sum_{\eta_{\Lambda_n}}\sum_{x \in \Lambda_n}\sum_{i\neq \eta_x}
        \left[\Gamma_n^\nu(x,\eta_x,\eta^{x,i})-\Gamma_n^\nu(x,i,\eta) \right]
        \log\left(\frac{\mu(\eta_{\Lambda_n})}{\mu(\eta_{\Lambda_n}^{x,i})}\right).
    \end{align*}
    Now we add and substract some terms to obtain
    \begin{align*}
    2 g^n_{\mathscr{L}}(\nu|\mu) 
    =
    -&\sum_{\eta_{\Lambda_n}}\sum_{x \in \Lambda_n}\sum_{i\neq \eta_x}
    \left[\Gamma_n^\nu(x,\eta_x,\eta^{x,i})-\Gamma_n^\nu(x,i,\eta) \right]
    \log\left(\frac{\Gamma_n^\nu(x, \eta_x, \eta^{x,i})}{\Gamma_n^\nu(x,i,\eta)}\right)
    \\\
    +
    &\sum_{\eta_{\Lambda_n}}\sum_{x \in \Lambda_n}\sum_{i\neq \eta_x}
    \left[\Gamma_n^\nu(x,\eta_x,\eta^{x,i})-\Gamma_n^\nu(x,i,\eta) \right]
    \\\
    &\quad \times \left\{
        \log\left(\frac{\nu(\eta_{\Lambda_n})}{\nu(\eta_{\Lambda_n}^{x,i})}\right)
        -
        \log\left(\frac{\mu(\eta_{\Lambda_n})}{\mu(\eta_{\Lambda_n}^{x,i})}\right)
        +
        \log\left(\frac{\Gamma_n^\nu(x, \eta_x, \eta^{x,i})}{\Gamma_n^\nu(x,i,\eta)}\right)
    \right\}
    \\\
    =
    -&\sum_{\eta_{\Lambda_n}}\sum_{x \in \Lambda_n}\sum_{i\neq \eta_x}
    \left[\Gamma_n^\nu(x,\eta_x,\eta^{x,i})-\Gamma_n^\nu(x,i,\eta) \right]
    \log\left(\frac{\Gamma_n^\nu(x, \eta_x, \eta^{x,i})}{\Gamma_n^\nu(x,i,\eta)}\right)
    \\\
    +
    &\sum_{\eta_{\Lambda_n}}\sum_{x \in \Lambda_n}\sum_{i\neq \eta_x}
    \left[\Gamma_n^\nu(x,\eta_x,\eta^{x,i})-\Gamma_n^\nu(x,i,\eta) \right]
    \\\
    &\quad 
    \times \left\{
        \log\left(\frac{\nu(\eta_{\Lambda_n})}{{\Gamma_n^\nu(x,i,\eta)}}\right)
        -
        \log\left(\frac{\nu(\eta^{x,i}_{\Lambda_n})}{\Gamma_n^\nu(x,\eta_x, \eta^{x,i})}\right)
        -
        \log\left(\frac{\mu(\eta_{\Lambda_n})}{\mu(\eta_{\Lambda_n}^{x,i})}\right)
    \right\},
\end{align*}
which yields the claim. 
\end{proof}
\subsection{A quantitative differentiation lemma}

To control the terms in the second sum in Lemma~\ref{lemma:finite-volume-relative-entropy-loss}, one can consider the following technical helper that tells us that the quotients of $\mu$ are actually approximating the conditional marginals of $\mu$ which are given by $\gamma$. Indeed, by using Lemma~\ref{differentiation-lemma} one can easily show the following convergence. 

\begin{lemma}\label{lemma:uniform-convergence-conditional-probabilities}
    Let $\mu \in \mathscr{G}(\gamma)$ and assume that the specification $\gamma$ is quasilocal and non-null. 
    Let $x \in \Z^d$ and fix $i\in \Omega_0$. Then, the following convergence holds uniform in $\eta \in \Omega$
    \begin{align*}
        \frac{\mu(\eta_{\Lambda_n})}{\mu(\eta^{x,i}_{\Lambda_n})} 
        \to 
        \frac{\gamma_x(\eta_x| \eta_{x^c})}{\gamma_x(i| \eta_{x^c})} 
        \quad 
        \text{as } n \to \infty.  
    \end{align*}
\end{lemma}

However, the above result does not give us any quantitative control over the speed of convergence. Therefore we have to use a tool that gives us a more precise result than Lemma \ref{differentiation-lemma}. 
For this, recall the notation
\begin{align*}
    \delta_x f  = \sup_{\eta \in \Omega, i \in \Omega_0}\abs{f(\eta^{x,i})-f(\eta)}.
\end{align*}

Now, as a first step towards getting a quantitative bound on the error term in Lemma~\ref{lemma:uniform-convergence-conditional-probabilities}, we obtain a quantitative version of Lemma~\ref{differentiation-lemma} in terms of the oscillations of $f$. 

\begin{lemma}\label{lemma:quantitative-differentiation}
    Let $\nu$ be a probability measure such that $\nu(\eta_{\Lambda})>0$ for all $\eta_{\Lambda}$. Then, for any function $f:\Omega \to \R$ with the property 
    \begin{align*}
        \sum_{x \in \Z^d}\delta_x f < \infty,
    \end{align*}
    the following uniform error estimate holds for all $\eta \in \Omega$
    \begin{align*}
        \abs{\frac{1}{\nu(\eta_\Lambda)}\int_{[\eta_\Lambda]}f(\omega)\nu(d\omega) - f(\eta)}
        \leq 
        \sum_{x \notin \Lambda}\delta_x f. 
    \end{align*}
\end{lemma}

\begin{proof}
    Fix $\eta \in \Omega$ and $\Lambda \Subset \Z^d$. Then we can fix an enumeration of the vertices in $\Lambda^c$ and write 
    $$[n] = \Lambda \cup \{x_1, \dots, x_n\}.$$ 
    By a telescope sum we see
    \begin{align*}
        f(\eta) - f(\omega)
        = \sum_{n=1}^\infty\left(f(\eta_{[n]}\omega_{[n]^c})-f(\eta_{[n-1]}\omega_{[n-1]^c})\right)
        \leq 
        \sum_{x \notin \Lambda}\delta_x f. 
    \end{align*}
    The claim now follows via integration. 
\end{proof}

Now let us apply this result to obtain a quantitative version of the convergence of conditional probabilities. Here we will use the short-hand notation
\begin{align*}
    \delta_y \gamma_x(\cdot) := \sup_{\omega \in \Omega, i \in \Omega_0}\abs{\gamma_x(\omega_x \lvert \omega_{x^c}^{y,i}) - \gamma_x(\omega_x \lvert \omega_{x^c})}, \quad x\neq y \in \Z^d. 
\end{align*}
This quantity tells us how much the conditional distribution of the particle at site $x$ depends on the state of the particle at a different site $y$. 
\begin{lemma}
     Let $\mu \in \mathscr{G}(\gamma)$ and assume that the specification $\gamma$ is quasilocal and non-null with constant $\delta >0$. Let $x \in \Z^d$ and fix $i \in \Omega_0$. Then, for all $\eta \in \Omega$ and $\Lambda \Subset \Z^d$, it holds that 
    \begin{align}
        \abs{\frac{\mu(\eta_{\Lambda})}{\mu(\eta_{\Lambda}^{x,i})}-\frac{\gamma_x(\eta_x \lvert \eta_{x^c})}{\gamma_x(i\lvert \eta_{x^c})}} \leq
        \frac{2}{\delta^2}\sum_{y \notin \Lambda}\delta_y \gamma_x(\cdot). 
    \end{align}
\end{lemma}

\begin{proof}
   As a first step, note that we can write
   \begin{align*}
       \frac{\mu(\eta_{\Lambda})}{\mu(\eta_{\Lambda_n}^{x,i})} = \frac{\mu(\eta_x \lvert \eta_{\Lambda \setminus \{x\}})}{\mu(i \lvert \eta_{\Lambda_n \setminus \{x\}})}.
   \end{align*}
   We first show uniform error bounds for the denominator and the numerator. For this, observe that the DLR equations imply 
   \begin{align*}
       \mu(\eta_x \lvert \eta_{\Lambda \setminus \{x\}}) = \frac{1}{\mu( \eta_{\Lambda \setminus \{x\}})}\int_{ [\eta_{\Lambda \setminus \{x\}}]}\gamma_x(\eta_x \lvert \omega_{x^c})\mu(d\omega). 
   \end{align*}
   To this  we can now apply the quantitative differentiation lemma to obtain 
   \begin{align*}
       \abs{ \mu(\eta_x \lvert \eta_{\Lambda \setminus  \{x\}}) - \gamma_x(\eta_x \lvert \eta_{x^c})} \leq \sum_{y \notin \Lambda}\delta_y \gamma_x(\cdot). 
   \end{align*}
   Analogously we obtain 
   \begin{align*}
       \abs{ \mu(i\lvert \eta_{\Lambda \setminus  \{x\}}) - \gamma_x(i \lvert \eta_{x^c})} \leq \sum_{y \notin \Lambda}\delta_y \gamma_x(\cdot). 
   \end{align*}
   Now we can use the simple algebraic rule 
   \begin{align*}
       ad - bc =  \frac{1}{2}[(a-b)(c+d) - (a+b)(c-d)]
   \end{align*}
   in conjunction with the non-nullness of $\gamma$, and hence $\mu$, to obtain the inequality
   \begin{align*}
       \abs{\frac{\mu(\eta_{\Lambda})}{\mu(\eta_{\Lambda_n}^{x,i})} - \frac{\gamma_x(\eta_x \lvert \eta_{x^c})}{\gamma_x(i\lvert \eta_{x^c})}}
       \leq 
       \frac{2}{\delta^2}\sum_{y \notin \Lambda}\delta_y \gamma_x(\cdot), 
   \end{align*}
   as desired.
\end{proof}
\subsection{The zero-loss equation}
The previously derived representation of the relative entropy loss in finite boxes $\Lambda_n$ directly implies the following equation in case the relative entropy loss vanishes. 
\begin{lemma}\label{lemma:zero-loss-equation}
    Let $n \in \N$ and $\nu \in \calM_1(\Omega)$ be such that there exists $T>0$ with 
    \begin{enumerate}[i.]
        \item $\nu P_s(\eta_{\Lambda_n})>0$ for all $\eta_{\Lambda_n}$ and $s \in [0,T]$  
        \item and $\int_0^T g^n_\mathscr{L}(\nu P_s  \lvert \mu) = 0$.
    \end{enumerate} 
    Then 
    \begin{align*}
        &\int_0^T\sum_{\eta_{\Lambda_n}}\sum_{x \in \Lambda_n}\sum_{i\neq \eta_x}
        \left[\Gamma_n^{\nu_s}(x,\eta_x,\eta^{x,i})-\Gamma_n^{\nu_s}(x,i,\eta) \right]
        \log\left(\frac{\Gamma_n^{\nu_s}(x, \eta_x, \eta^{x,i})}{\Gamma_n^{\nu_s}(x,i,\eta)}\right)ds
        \\\
        =
        \int_0^T
        &\sum_{\eta_{\Lambda_n}}\sum_{x \in \Lambda_n}\sum_{i\neq \eta_x}
        \left[\Gamma_n^{\nu_s}(x,\eta_x,\eta^{x,i})-\Gamma_n^{\nu_s}(x,i,\eta) \right]
        \\\
        &\quad 
        \times \left\{
        \log\left(\frac{\nu_s(\eta_{\Lambda_n})}{{\Gamma_n^{\nu_s}(x,i,\eta)}}\right)
        -
        \log\left(\frac{\nu_s(\eta^{x,i}_{\Lambda_n})}{\Gamma_n^{\nu_s}(x,\eta_x, \eta^{x,i})}\right)
        -
        \log\left(\frac{\mu(\eta_{\Lambda_n})}{\mu(\eta_{\Lambda_n}^{x,i})}\right)
        \right\}ds.  
\end{align*}
\end{lemma}

We now want to estimate the terms appearing in this equation with the final goal to show that every term on the right-hand side actually vanishes. 
For this, let us first introduce some more notation, 
\begin{align*}
    \alpha_n(x,\nu) &= \sum_{\eta_{\Lambda_n}}\sum_{i \neq \eta_x}\left[\Gamma_n^\nu(x,\eta_x,\eta^{x,i})-\Gamma_n^\nu(x,i,\eta) \right]
    \log\left(\frac{\Gamma_n^\nu(x, \eta_x, \eta^{x,i})}{\Gamma_n^\nu(x,i,\eta)}\right), 
    \\\
    \beta_n(x,\nu) &= \sum_{\eta_{\Lambda_n}}\sum_{i \neq \eta_x}\lvert\Gamma_n^\nu(x,\eta_x,\eta^{x,i})-\Gamma_n^\nu(x,i,\eta) \rvert,
    \\\
    \rho_n(x) &= \sum_{y \notin \Lambda_n}\left(\delta_y \gamma_x(\cdot) + \delta_y c_x(\cdot)\right). 
\end{align*}

\begin{lemma}\label{lemma:time-averaged-upperbound-alpha}
    Let $n \in \N$ and $\nu \in \calM_1(\Omega)$ be such that there exists $T>0$ with 
    \begin{enumerate}[i.]
        \item $\nu P_s(\eta_{\Lambda_n})>0$ for all $\eta_{\Lambda_n}$ and $s \in [0,T]$  
        \item and $\int_0^T g^n_\mathscr{L}(\nu P_s  \lvert \mu) = 0$.
    \end{enumerate} 
    Then there exists a constant $C > 0$ that does not depend on $n$ and $\nu$ such that 
    \begin{align*}
        \int_0^T\sum_{x \in \Lambda_n}\alpha_n(x,\nu_s) ds\leq C \int_0^T\sum_{x\in \Lambda_n}\beta_n(x,\nu_s)\rho_n(x)ds. 
    \end{align*}
\end{lemma}

\begin{proof}
By Lemma~\ref{lemma:zero-loss-equation} it suffices to show that, for all $s \in [0,T]$, $x \in \Lambda_n$,  and $ \eta_{\Lambda_n}$, we have
    \begin{align*}
           \left\lvert 
        \log\left(\frac{\nu_s(\eta_{\Lambda_n})}{{\Gamma_n^{\nu_s}(x,i,\eta)}}\right)
        -
        \log\left(\frac{\nu_s(\eta^{x,i}_{\Lambda_n})}{\Gamma_n^{\nu_s}(x,\eta_x, \eta^{x,i})}\right)
        -
        \log\left(\frac{\mu(\eta_{\Lambda_n})}{\mu(\eta_{\Lambda_n}^{x,i})}\right)
        \right\rvert 
        \leq 
        C
        \rho_n(x)
    \end{align*}
    for some $C>0$. 
    To do this, we introduce a reference configuration $\Bar{\eta}$ defined by 
    \begin{align*}
        \overline{\eta} 
        =
        \begin{cases}
        \eta_y, \quad &\text{if $y \in \Lambda_n$}, \\\
        1, &\text{otherwise}. 
        \end{cases}
    \end{align*}
    and add and subtract terms of the form 
    \begin{align*}
        \log\left(c_x(\overline{\eta},i)\gamma_x(\eta_x \lvert \overline{\eta}_{x^c}) \right).
    \end{align*}
    By detailed balance we have
    \begin{align*}
        c_x(\overline{\eta},i)\gamma_x(\eta_x \lvert \overline{\eta}_{x^c})
        =
        c_x(\overline{\eta}^{x,i},\eta_x)\gamma_x(i \lvert \overline{\eta}_{x^c}). 
    \end{align*}
    So we actually just need to add (or subtract)
    \begin{align*}
        0 = \log\left(\frac{c_x(\overline{\eta},i)\gamma_x(\eta_x \lvert \overline{\eta}_{x^c})}{c_x(\overline{\eta}^{x,i},\eta_x)\gamma_x(i \lvert \overline{\eta}_{x^c})}\right).
    \end{align*}
    For every term in the sum above this gives us 
    \begin{align*}
        &\left\lvert
        \log\left(\frac{\nu(\eta_{\Lambda_n})}{{\Gamma_n^\nu(x,i,\eta)}}\right)
        -
        \log\left(\frac{\nu(\eta^{x,i}_{\Lambda_n})}{\Gamma_n^\nu(x,\eta_x, \eta^{x,i})}\right)
        -
        \log\left(\frac{\mu(\eta_{\Lambda_n})}{\mu(\eta_{\Lambda_n}^{x,i})}\right)
        +
        \log\left(\frac{c_x(\overline{\eta},i)\gamma_x(\eta_x \lvert \overline{\eta}_{x^c})}{c_x(\overline{\eta}^{x,i},\eta_x)\gamma_x(i \lvert \overline{\eta}_{x^c})}\right)
        \right\rvert 
        \\\
        =
        &\left\lvert 
        \log\left(\frac{\nu(\eta_{\Lambda_n})c_x(\overline{\eta},i)}{{\Gamma_n^\nu(x,i,\eta)}} \right)
        -
\log\left(\frac{\nu(\eta^{x,i}_{\Lambda_n})c_x(\overline{\eta}^{x,i},\eta_x)}{\Gamma_n^\nu(x,\eta_x, \eta^{x,i})}\right)
        -
        \log\left(\frac{\mu(\eta_{\Lambda_n})}{\mu(\eta_{\Lambda_n}^{x,i})}\frac{\gamma_x(i\lvert \overline{\eta}_{x^c})}{\gamma_x(\eta_x \lvert \overline{\eta}_{x^c})}\right)
        \right\rvert.
    \end{align*}
    By assumption $\mathbf{(R3)}$ and $\mathbf{(S2)}$ we know that all the terms in the logarithms are bounded away from $0$, so we can use the Lipschitz continuity of $\log(\cdot)$ away from $0$, apply  Lemma~\ref{lemma:quantitative-differentiation} to the first two terms and Lemma~\ref{lemma:uniform-convergence-conditional-probabilities} to the third term to obtain the estimate 
    \begin{align*}
        &\left\lvert 
        \log\left(\frac{\nu_s(\eta_{\Lambda_n})c_x(\overline{\eta},i)}{{\Gamma_n^{\nu_s}(x,i,\eta)}} \right)
        -\log\left(\frac{\nu_s(\eta^{x,i}_{\Lambda_n})c_x(\overline{\eta}^{x,i},\eta_x)}{\Gamma_n^{\nu_s}(x,\eta_x, \eta^{x,i})}\right)
        -
        \log\left(\frac{\mu(\eta_{\Lambda_n})}{\mu(\eta_{\Lambda_n}^{x,i})}\frac{\gamma_x(i\lvert \overline{\eta}_{x^c})}{\gamma_x(\eta_x \lvert \overline{\eta}_{x^c})}\right)
        \right\rvert
        \\\
        \leq 
        &L\Bigg(
        \abs{c_x(\Bar{\eta},i)-\frac{\Gamma_n^{\nu_s}(x,i,\eta)}{\nu_s(\eta_{\Lambda_n})}}
        +
        \abs{c_x(\Bar{\eta}^{x,i},\eta_x) - \frac{\Gamma_n^{\nu_s}(x,\eta_x, \eta^{x,i})}{\nu_s(\eta^{x,i}_{\Lambda_n})}}
        +
        \abs{\frac{\mu(\eta_{\Lambda_n})}{\mu(\eta_{\Lambda_n}^{x,i})} - \frac{\gamma_x(\eta_x \lvert \bar{\eta}_{x^c})}{\gamma_x(i\lvert \Bar{\eta}_{x^c})}}\Bigg)
        \\\
        \leq 
        &\frac{2L}{\delta^2} \sum_{y \notin \Lambda_n} \left( \delta_y c_x(\cdot) + \delta_y \gamma_x(\cdot) \right), 
    \end{align*}
    where $\delta >0$ is the non-nullness constant of the specification $\gamma$.  
\end{proof}

As a next step,  we will show that, for fixed $\nu$, the quantity $\alpha_n(x,\nu)$ is non-decreasing in $n$. Also note that each summand in the definition of $\alpha_n(x,\nu)$ is actually non-negative.  
\begin{lemma}\label{lemma:monotonicity-alpha}
    For any $\nu \in \calM_1(\Omega)$ such that $\nu(\eta_{\Lambda})>0$ for all $\Lambda \Subset \Z^d$ and $\eta_{\Lambda} \in \Omega_\Lambda$, it holds that for all $n \in \N$ and $x \in \Lambda_n \subset  \Lambda_{n+1}$ 
    \begin{align*}
        0 \leq \alpha_n(x,\nu) \leq \alpha_{n+1}(x,\nu). 
    \end{align*}
\end{lemma}
\begin{proof}
    Define the function 
    \begin{align*}
        \Phi(u,v) = (u-v)\log\left(\frac{u}{v}\right), \quad u,v >0. 
    \end{align*}
    Then $\Phi$ is convex and homogeneous of degree one, i.e., $\Phi(\lambda u, \lambda v) = \lambda \Phi(u,v)$ for all $\lambda > 0$. This implies that it is subadditive. Indeed, for all $u_1, u_2, v_1, v_2>0$ we have that 
    \begin{align*}
        \Phi(u_1 + u_2, v_1 + v_2) 
        &= 
        2\Phi\left( \frac{1}{2}u_1+\frac{1}{2}u_2, \frac{1}{2}v_1 + \frac{1}{2}v_2\right)
        \\\
        &\leq 
        2\left[ \frac{1}{2}\Phi(u_1, v_1) + \frac{1}{2}\Phi(u_2, v_2)\right]
        =
        \Phi(u_1,v_2) + \Phi(u_2, v_2). 
    \end{align*}
    We can rewrite $\alpha_n(x,\nu)$ and $\alpha_{n+1}(x,\nu)$ in terms of $\Phi$ as 
    \begin{align*}
        \alpha_n(x,\nu) &= \sum_{\eta_{\Lambda_n}}\sum_{i \neq \eta_x}\Phi\left(\Gamma_n^\nu(x, i, \eta_{\Lambda_n}), \Gamma_n^\nu(x, \eta_x, \eta^{x,i}_{\Lambda_n})\right),\\\
        \alpha_{n+1}(x,\nu) &= \sum_{\eta_{\Lambda_{n+1}}}\sum_{i \neq \eta_x}\Phi\left(\Gamma^\nu_{n+1}(x, i, \eta_{\Lambda_n}), \Gamma^\nu_{n+1}(x, \eta_x, \eta^{x,i}_{\Lambda_{n+1}})\right). 
    \end{align*}
    Since we have 
     \begin{align*}
         \Gamma_n^\nu(x,i,\eta_{\Lambda_n}) = \sum_{\xi_{\Lambda_{n+1}}: \xi_{\Lambda_{n}} = \eta_{\Lambda_n}}\Gamma^\nu_{n+1}(x,i,\xi_{\Lambda_{n+1}}),
     \end{align*}
     the claim follows from the subadditivity of $\Phi$. 
\end{proof}
Let us briefly pause here and discuss what we have shown so far and where we are headed. We want to show that all of the terms $\alpha_n(x)$ vanish and in Lemma~\ref{lemma:time-averaged-upperbound-alpha} we have established the inequality 
\begin{align*}
    0 \leq \int_0^T\sum_{x \in \Lambda_n}\alpha_n(x,\nu_s) ds\leq C \int_0^T\sum_{x\in \Lambda_n}\beta_n(x,\nu_s)\rho_n(x)ds. 
\end{align*}
Now if we were able to show that we can control the terms $\beta_n(x,\cdot)$ in terms of $\alpha_n(x,\cdot)$, then a sufficiently fast decay of $\rho_n(x)$ should allow us to conclude that the $\alpha_n(x,\cdot)$ vanish. But this decay is implied by our assumptions.  
Therefore, our strategy will be the following. We first establish a pointwise estimate for $\beta_n(x,\cdot)$ in terms of $\alpha_n(x,\cdot)$ and then put everything together at the end of the section to prove Proposition~\ref{proposition:time-averaged-holley-stroock}. 

\begin{lemma}\label{lemma:upper-bound-beta}
    Let $\nu \in \calM_1(\Omega)$. Then for all $n \in \N$ and $x \in \Lambda_n$ it holds that 
    \begin{align*}
        \beta_n(x,\nu)^2 \leq |\Omega_0| \cdot \sup_{\omega, i}\abs{c_x(\omega,i)}\alpha_n(x,\nu). 
    \end{align*}

\end{lemma}

\begin{proof}
    Here we can use the symmetry and subadditivity of $\Phi$ to show that 
    \begin{align*}
        \alpha_n(x,\nu)
        =
        \sum_{\eta_{\Lambda_n}}\sum_{i\neq \eta_x}\Phi\left(\Gamma_n^\nu(x,i,\eta_{\Lambda_n}), \Gamma(x,\eta_x, \eta_{\Lambda_n}^{x,i})\right)
        \geq 
        \Phi(M,m),
    \end{align*}
    where we use the notation
    \begin{align*}
        M &= \sum_{\eta_{\Lambda_n}}\sum_{i\neq \eta_x}\max\left\{\Gamma_n^\nu(x,i, \eta_{\Lambda_n}), \Gamma_n^\nu(x,\eta_x, \eta^{x,i}_{\Lambda_n})\right\},
        \\\
        m &= \sum_{\eta_{\Lambda_n}}\sum_{i\neq \eta_x}\min\left\{\Gamma_n^\nu(x,i, \eta_{\Lambda_n}), \Gamma_n^\nu(x,\eta_x, \eta^{x,i}_{\Lambda_n})\right\}. 
    \end{align*}
    Since we have 
    \begin{align*}
        \beta_n(x,\nu) = M - m
    \end{align*}
    and the trivial bound 
    \begin{align*}
        M \leq |\Omega_0| \sup_{\omega, i} c_x(\omega,i), 
    \end{align*}
    the claimed inequality follows from the fact that
    \begin{align*}
        u - v \leq u \log \frac{u}{v}, \quad 0 < v \leq u. 
    \end{align*}
    Indeed, by the above calculations we have 
    \begin{align*}
        \alpha_n(x,\nu) \cdot M \geq (M-m)\log\left(\frac{M}{m}\right) \cdot M 
        \geq (M-m)^2 = \beta_n(x,\nu)^2. 
    \end{align*}
    Using the previously derived upper bound on $M$ now yields the claim. 
\end{proof}


With all of these rather technical estimates in place, we are finally ready to prove Proposition~\ref{proposition:time-averaged-holley-stroock}. 

\begin{proof}[Proof of Proposition \ref{proposition:time-averaged-holley-stroock}]
    Recall that our goal is to show that $\alpha_n(x,\cdot) \equiv 0$ for every $n \in \N$ and $x\in \Lambda_n$. As a first step we will show that it suffices to prove that 
\begin{align*}
    \mathbf{C}_1 := \sup_{x \in \Z^d} \sum_{n=1}^\infty \rho_n(x) < \infty\quad\text{ and }\quad \mathbf{C}_2 := \sup_{n \in \N} \frac{1}{n^{d-1}}\sum_{x \in \Lambda_n} \rho_n(x) < \infty. 
\end{align*}
    Indeed, if this is the case, we can first combine Lemma \ref{lemma:time-averaged-upperbound-alpha} and the pointwise estimates from Lemma \ref{lemma:upper-bound-beta} to obtain
    \begin{align}\label{ineq:initial-bound}
        \int_0^T \sum_{x \in \Lambda_n}\alpha_n(x,\nu_s)ds 
        \leq 
        C' \int_0^T \sum_{x \in \Lambda_n}\beta_n(x,\nu_s)\rho_n(x)ds
        \leq 
        C
        \int_0^T \sum_{x \in \Lambda_n}\sqrt{\alpha_n(x,\nu_s)}\rho_n(x)ds. 
    \end{align}
    By Lemma \ref{lemma:monotonicity-alpha} we have the following pointwise estimate for all $s\in [0,T]$
    \begin{align}\label{ineq:pointwise-alpha-sum-lower-bound} 
    \sum_{x \in \Lambda_n}\alpha_n(x,\nu_s) 
    &\geq \mathbf{C}_1^{-1}\sum_{x \in \Lambda_n}\alpha_n(x,\nu_s)\sum_{k=1}^n \rho_k(x)
   \geq \mathbf{C}_1^{-1}\sum_{k=1}^n\sum_{x \in \Lambda_k} \alpha_k(x,\nu_s)\rho_k(x). 
    \end{align}
    Since the coefficients $\rho_n(x)$ do not depend on $\nu_s$, we can pull them out of the integrals and define for $k \in \N$
    \begin{align*}
        \delta_k := 
        \sum_{x \in \Lambda_k}\rho_k(x) \int_0^T \alpha_k(x,\nu_s)ds. 
    \end{align*}
    Note that by definition of $\alpha$ and $\rho$ we have $\delta_k \geq 0$ for all $k \in \N$. By combining \eqref{ineq:pointwise-alpha-sum-lower-bound} and \eqref{ineq:initial-bound} with the Cauchy--Schwarz inequality for sums we obtain 
    \begin{align*}
        \left[\sum_{k=1}^n \delta_k\right]^2 \leq C^2 \sum_{x \in \Lambda_n}\rho_n(x)\sum_{x \in \Lambda_n}\rho_n(x)\left(\int_0^T \sqrt{\alpha_n(x,\nu_s)}ds\right)^2. 
    \end{align*}
    Another application of the Cauchy--Schwarz inequality to the integrals on the right-hand side yields 
    \begin{align*}
        \left(\int_0^T \sqrt{\alpha_n(x,\nu_s)}ds\right)^2 
        \leq 
        T \int_0^T \alpha_n(x,\nu_s)ds. 
    \end{align*}
    So we finally obtain 
    \begin{align}\label{ineq:final-estimate}
        \left[\sum_{k=1}^n\delta_k\right]^2 \leq T \mathbf{C}_2 \mathbf{C}_1^2 C^2 \delta_n n^{d-1} =: \mathbf{C}\delta_n n^{d-1}. 
    \end{align}
    If there were an index $n_0 \in \N$ such that $\delta_{n_0}>0$, then for all $n>n_0$ it would hold that 
    \begin{align*}
        \frac{1}{n^{d-1}} \leq \mathbf{C}\left[\frac{1}{\sum_{k=1}^{n-1}\delta_k}-\frac{1}{\sum_{k=1}^n \delta_k}\right]. 
    \end{align*}
    By a standard telescoping argument and monotonicity, the series over the terms on the right-hand side converges, which leads to a contradiction for $d \in \{1,2\}$. Therefore, we must have $\delta_n = 0$ for all $n \in \N$ and hence by continuity $\alpha_n(x,\cdot) \equiv 0$ for all $n \in \N$ and $x \in \Lambda_n$. So it remains to show that $\mathbf{C}_1$ and $\mathbf{C}_2$ are actually finite. 
    \medskip 

\noindent
\textit{Ad $\mathbf{C}_1$}: For fixed $x\in\Z^d$ we have 
\begin{align*}
    \sum_{n=1}^\infty \rho_n(x)
    &=
    \sum_{n=1}^\infty \sum_{y \notin \Lambda_n}(\delta_y \gamma_x(\cdot) + \delta_y c_x(\cdot))
    \\\
    &=
    \sum_{y \in \Z^d}(\delta_y \gamma_x(\cdot) + \delta_y c_x(\cdot))\abs{\{n \in \N: \ x \in \Lambda_n, y \notin \Lambda_n\}}
    \\\
    &\leq 
    \sum_{y \in \Z^d}(\delta_y \gamma_x(\cdot) + \delta_y c_x(\cdot)) \abs{x-y}.
\end{align*}
Now assumptions $\mathbf{(R4)}$ and $\mathbf{(S3)}$ yield a uniform in $x$ upper bound on this quantity. 

\medskip 
\noindent 
\textit{Ad $\mathbf{C}_2$}: Here we have for fixed $n \in  \N$
\begin{align*}
    \sum_{x \in  \Lambda_n}\sum_{y \notin \Lambda_n}(\delta_y \gamma_x(\cdot) + \delta_y c_x(\cdot))
    &\leq
    \sum_{v \in \Z^d}\sum_{x \in \Lambda_n: \ x+v \notin \Lambda_n}(\delta_{x+v}\gamma_x(\cdot) + \delta_{x+v}c_x(\cdot))
    \\\
    &\leq 
    d(2n+1)^{d-1} \sum_{v \in \Z^d}\abs{v}\sup_{x \in \Z^d}(\delta_{x+v}\gamma_x(\cdot) + \delta_{x+v}c_x(\cdot)).
\end{align*}
This can be bounded from above, independent of $n$, by assumptions $\mathbf{(R4)}$ and $\mathbf{(S3)}$. 
\end{proof}

\begin{remark}
    Only the very last step of the proof of Proposition \ref{proposition:time-averaged-holley-stroock} depends on the dimension $d$, so all of the estimates up to this point, including \eqref{ineq:final-estimate}, hold in any dimension.
    Therefore, let us take another look at this key estimate before we move on.
    If one had a uniform lower bound on $\alpha_{l}(x)$ for all $x \in \Z^d$ and some $l \in \N$, then, assuming that the coefficients $\rho$ are non-trivial, we would have $\delta_n \sim n^{d-1}$ and hence for sufficiently large $n$
    \begin{align*}
        \left[\sum_{k=1}^n \delta_k\right]^2 \sim n^{2d} 
        \quad
        \text{and}
        \quad 
        \delta_n n^{d-1} \sim n^{2d-2}. 
    \end{align*}
    So no matter what the constants on the left-hand side of \eqref{ineq:final-estimate} are, this estimate directly gives us a contradiction for sufficiently large $n$. Now if $\nu \in \calM_1(\Omega)$ is stationary but not reversible, then we get this uniform control over $\alpha_l(x)$ under the additional assumption of translation-invariance. 
    However, without this additional assumption it is not clear, how to show that the right-hand side actually grows sufficiently fast to obtain this contradiction.  
\end{remark}

\section{Proof of the positive-mass property}\label{section:positive-mass}

We actually show the following more general result which implies Proposition~\ref{proposition:positive-mass} as a special case.

\begin{proposition}
Consider an interacting particle system with single-site updates with generator given by 
\begin{align*}
    Lf(\eta) = \sum_{i  \in \Z^d}\sum_{\xi_i = 1, \dots, q}c_i(\eta, \xi_i)\left[f(\xi_i \eta_{i^c}) - f(\eta)\right],
\end{align*}
where the rates satisfy the assumptions $(\mathbf{L1})-(\mathbf{L2})$.
We further assume the following.
\begin{enumerate}[\bfseries \upshape (R1')]
    \item For $L$, reachability is independent of the boundary conditions, i.e., whenever $c_x(\eta, \xi_x) >0$ we also have $c_x(\sigma, \xi_x) >0$ for all $\sigma$ with $\sigma_x = \eta_x$. 
    In this case we say that $\xi_x$ is reachable from $\eta_x$ and write $d_x(\eta_x, \xi_x)$ for the indicator of this event. 
    \item $L$ is single-site irreducible, i.e., the Markov chain on the single state space with rates given by $d$ is irreducible. 
    \item The minimal transition rate is strictly positive, i.e.,
    \begin{align*}
        \hat{c} = \inf_{i, \omega, \xi_i: \ c_i(\omega, \xi_i) >0}c_i(\omega, \xi_i) > 0.
    \end{align*}
 \end{enumerate}
 Then, for all $\tau >0$ and $\Lambda \Subset \Z^d$, there exists a constant $C(\tau, \Lambda) > 0$ such that for any initial distribution $\nu \in \calM_1(\Omega)$ and any time $t \in [\tau, \infty)$ we have 
 \begin{align*}
     \forall \eta \in \Omega: \quad \nu_t(\eta_\Lambda) \geq C(\tau, \Lambda). 
 \end{align*}
 In particular, for all subsequential limits $\nu^* = \lim_{n \to \infty} \nu_{t_n}$ with $t_n \uparrow \infty$ we have 
 \begin{align*}
     \forall \eta \in \Omega \ \forall \Lambda \Subset \Z^d: \quad \nu^*(\eta_\Lambda) > 0. 
 \end{align*}
\end{proposition}

\begin{proof}
We will show that there exists $C(\tau, \Lambda) > 0$ such that for all initial infinite-volume configurations $\omega$, there is a lower bound 
\begin{align*}
    \mathbb{Q}_\omega[\sigma_\Lambda(\tau) = \eta_\Lambda] \geq C(\tau, \Lambda). 
\end{align*}
This then implies for all $t \in [\tau, \infty)$ and $\nu \in \calM_1(\Omega)$
\begin{align*}
    \nu_t(\eta_\Lambda) = \int_\Omega \mathbb{Q}_\omega[\sigma_\Lambda(\tau) = \eta_\Lambda]\nu_{t - \tau}(d\omega) \geq C(\tau, \Lambda). 
\end{align*}
\medskip 
    
\noindent 
\textit{Step 0: } We will compare our dynamics to a second interacting particle system with generator $\hat{L}$, that can be seen as a finite-volume perturbation of $L$. More precisely, we consider a generator $\hat{L}$ with rates $\hat{c}$ that agree with the rates $c$, except inside of $\Lambda$, where all of the sites $i \in \Lambda$ will behave independently with transition rates 
$\hat{c}_i(\eta, \xi_i) = \hat{c} \cdot d_i(\eta_i, \xi_i)$.
We will denote the path measure on the space of $\Omega$-valued cádlág paths with respect to the original dynamics by $\mathbb{Q}$ and with respect to the perturbed dynamics by $\hat{\mathbb{Q}}^\Lambda$. 
\medskip 

\noindent 
\textit{Step 1:} By Lemma~\ref{lemma:girsanov-transformation} we have the following  Girsanov-type formula on the space of cádlág paths
\begin{align*}
    \frac{d \mathbb{Q}_\omega}{d\hat{\mathbb{Q}}^\Lambda_\omega} ( \sigma[0,\tau])
    =
    \exp\left(-\int_0^\tau \lambda(\sigma(s))ds + \sum_{s \in [0,\tau]: \sigma_{\Lambda}(s_-) \neq \sigma_{\Lambda}(s)}\sum_{i \in \Lambda}\log\left(\frac{c_i(\sigma(s_-), \sigma_i(s))}{\hat{c}_i(\sigma(s_-), \sigma_i(s))}\right)\right), 
\end{align*}
where 
\begin{align*}
    \lambda(\eta) :=\sum_{i \in \Lambda}\left(c_i(\eta) - \hat{c}_i(\eta)\right),
\end{align*}
and $c_i(\eta)$ respectively $\hat{c}_i(\eta)$ are the total rates at which we see a flip at site $i$ when we are currently in configuration $\eta$, i.e., 
\begin{align*}
        c_i(\eta) = \sum_{\xi_i \neq \eta_i}c_i(\eta, \xi_i) \quad \text{respectively} \quad \hat{c}_i(\eta) = \sum_{\xi_i \neq \eta_i}\hat{c}_i(\eta, \xi_i). 
\end{align*}
Let us introduce additional notation to refer to the two separate parts of the above Radon--Nikodym derivative 
\begin{align*}
    a(\sigma[0,\tau]) &= \exp\left(-\int_0^\tau \lambda(\sigma(s))ds\right),
    \\\
    A(\sigma[0,\tau]) &= \exp\left(\sum_{s \in [0,\tau]: \sigma_{\Lambda}(s_-) \neq \sigma_{\Lambda}(s)}\sum_{i \in \Lambda}\log\left(\frac{c_i(\sigma(s_-), \sigma_i(s))}{\hat{c}_i(\sigma(s_-), \sigma_i(s))}\right)\right). 
\end{align*}
\medskip

\noindent 
\textit{Step 2: } By the Girsanov-type formula we can rewrite the probability we want to bound as
\begin{align*}
        \mathbb{Q}_\omega(\sigma_\Lambda(\tau) = \eta_\Lambda(\tau))
        =
        \hat{\mathbb{Q}}^\Lambda_\omega\left(a(\sigma[0,\tau]) A(\sigma[0,\tau]) \mathbf{1}_{\{\sigma_\Lambda(\tau) = \eta_\Lambda\}}\right). 
\end{align*}
  So in order to obtain the claimed lower bound, we only need to lower bound the functionals $a$ and $A$.
    \medskip

    \noindent
    \textit{Step 3: Deterministic lower bound on $a$.} 
    Since the rates of $L$ and the finite-volume perturbation $\hat{L}$ are assumed to be bounded from above by some constant $\mathbf{c}$ (this follows from well-definedness via Liggetts criteria), we can bound the function $\lambda(\cdot)$ by 
    \begin{align*}
    - \mathbf{c}\abs{\Lambda} \leq \lambda(\sigma[0,\tau]) \leq \mathbf{c}\abs{\Lambda}.
    \end{align*}
    This translates into an upper and lower bound for $a(\sigma[0,\tau])$ via 
    \begin{align*}
    \kappa(\tau) := \exp(-\tau \mathbf{c}\abs{\Lambda}) \leq a([0,\tau]) \leq \exp(\tau \mathbf{c}\abs{\Lambda}). 
    \end{align*}
    This implies 
    \begin{align*}
        \hat{\mathbb{Q}}^\Lambda_\omega\left(a(\sigma[0,\tau]) A(\sigma[0,\tau]) \mathbf{1}_{\{\sigma_\Lambda(\tau) = \eta_\Lambda\}}\right)
        \geq 
        \kappa(\tau, \Lambda) 
        \hat{\mathbb{Q}}^\Lambda_\omega\left(A(\sigma[0,\tau]) \mathbf{1}_{\{\sigma_\Lambda(\tau) = \eta_\Lambda\}}\right). 
    \end{align*}
    \medskip 

    \noindent 
    \textit{Step 4: Probabilistic lower bounds for $A$.} 
    We can lower bound $A$ in terms of the total number of jumps inside $\Lambda$, i.e., 
    \begin{align*}
        A(\sigma[0,\tau]) \geq e^{-R N_\Lambda(\tau)},
    \end{align*}
    where $R>0$ is a constant that depends on the minimal positive and maximal transition rate of $L$ and for any volume $\Delta \Subset \Z^d$ we define the (almost-surely finite) random variable $N_\Delta(\tau)$ by
    \begin{align*}
        N_\Delta(\tau) := \abs{\{s \in [0,\tau]: \ \sigma(s_-) \neq \sigma(s)\}}.
    \end{align*}
    In this notation we obtain 
    \begin{align*}
        \hat{\mathbb{Q}}^\Lambda_\omega\left(A(\sigma[0,\tau]) \mathbf{1}_{\{\sigma_\Lambda(\tau) = \eta_\Lambda\}}\right)
        \geq 
        \hat{\mathbb{Q}}^\Lambda_\omega\left(e^{-RN_\Lambda(\tau)} \mathbf{1}_{\{\sigma_\Lambda(\tau) = \eta_\Lambda\}}\right).
    \end{align*}
    \medskip 

    \noindent 
    \textit{Step 5: Using the independence to factorize.}
    Now note that under $\hat{\mathbb{Q}}^\Lambda_\omega$ all spins inside of $\Lambda$ are independent, and we have $N_\Lambda = \sum_{i\in \Lambda}N_i$. So we should get a really nice factorization, more precisely 
    \begin{align*}
        \hat{\mathbb{Q}}^\Lambda_\omega\left(e^{-RN_\Lambda(\tau)} \mathbf{1}_{\{\sigma_\Lambda(\tau) = \eta_\Lambda\}}\right)
        = \prod_{i \in \Lambda}\hat{\mathbb{Q}}^\Lambda_\omega\left(e^{-RN_i(\tau)} \mathbf{1}_{\{\sigma_i(\tau) = \eta_i\}}\right).
    \end{align*}
    \medskip 

    \noindent 
    \textit{Step 6: Estimating the factors via Poisson tails.}
    The factors can now each be estimated separately. For every $i \in \Lambda$ we have
    \begin{align*}
        \hat{\mathbb{Q}}^\Lambda_\omega\left(e^{-RN_i(\tau)} \mathbf{1}_{\{\sigma_i(\tau) = \eta_i\}}\right)
        &\geq 
        e^{-Rm}\hat{\mathbb{Q}}^\Lambda_\omega\left(N_0(\tau) \leq m, \sigma_i(\tau) = \eta_i \right)
        \\\
        &\geq
        e^{-Rm} \left(\hat{\mathbb{Q}}^\Lambda_\omega\left(\sigma_i(\tau) = \eta_i \right) - \hat{\mathbb{Q}}^\Lambda_\omega\left(N_0(\tau) > m \right)\right).
    \end{align*}
    By irreducibility of the single-site dynamics, there is a strictly positive lower bound for $\hat{\mathbb{Q}}^\Lambda_\omega\left(\sigma_i(\tau) = \eta_i \right)$ that only depends on $\tau$ (and not on $\omega$, $i$ or $\eta_i$) and the second term is the tail of a Poisson random variable. Hence, we will need to choose $m$ sufficiently large to make the right-hand side positive. So let us choose such an $m$ and denote the thereby obtained lower bound by $\rho(\tau) >0$.
    
    By putting all of the steps above together the claimed lower bound follows. 
\end{proof}


\section{Towards a generalisation}\label{section:generalisations}
\subsection{Generalisation to synchronous multi-site updates}
While the smoothness assumptions on the rates and the specification seem to be quite natural, it would be nice to lift the restriction of only considering single-site updates. However, there one runs into the trouble that 
\begin{align*}
    g_n^\mathscr{L}(\nu \lvert \mu) 
    = 
    &\sum_{\Delta \cap \Lambda_n \neq \emptyset} \sum_{\xi_\Delta}\int_\Omega \nu(d\eta) c_\Delta(\eta,\xi_\Delta)\log \frac{\nu(\xi_{\Delta \cap \Lambda_n}\eta_{\Lambda_n \setminus \Delta}) \mu (\eta_{\Lambda_n})}{\mu(\xi_{\Delta \cap \Lambda_n}\eta_{\Lambda_n \setminus \Delta})\nu(\eta_{\Lambda_n})}
    \\\
    =
    &\sum_{\Delta \subset \Lambda_n} \sum_{\xi_\Delta}\int_\Omega \nu(d\eta) c_\Delta(\eta,\xi_\Delta)\log \frac{\nu(\xi_{\Delta \cap \Lambda_n}\eta_{\Lambda_n \setminus \Delta}) \mu (\eta_{\Lambda_n})}{\mu(\xi_{\Delta \cap \Lambda_n}\eta_{\Lambda_n \setminus \Delta})\nu(\eta_{\Lambda_n})}
    \\\
    +
    &\sum_{\Delta \nsubseteq \Lambda_n : \Delta \cap \Lambda_n \neq \emptyset}\sum_{\xi_\Delta}\int_\Omega \nu(d\eta) c_\Delta(\eta,\xi_\Delta)\log \frac{\nu(\xi_{\Delta \cap \Lambda_n}\eta_{\Lambda_n \setminus \Delta}) \mu (\eta_{\Lambda_n})}{\mu(\xi_{\Delta \cap \Lambda_n}\eta_{\Lambda_n \setminus \Delta})\nu(\eta_{\Lambda_n})}.
\end{align*}
We can only rewrite the first sum as in Lemma \ref{lemma:finite-volume-relative-entropy-loss}, but this doesn't work for the second sum, because there are not enough terms to perform the change of variables.
One could now try to just bound the second term and thereby obtain Lemma \ref{lemma:time-averaged-upperbound-alpha} with an additional additive error term on the right-hand side. One can show that this error term is of boundary order, see \cite[Lemma 3.10]{jahnel_dynamical_2022} for a proof in a more general setting, but naively carrying this term through the rest of the proof makes it impossible to apply the summability argument at the end of the proof of Proposition~\ref{proposition:time-averaged-holley-stroock}. 

This is a bit strange, because at least heuristically one could say that the Holley--Stroock argument works in one and two dimensions, because the boundary contributions cannot weigh up against the bulk contribution and the error term is also of boundary order, but unfortunately we have not been able to make this (potentially misguided) intuition rigorous. 

\subsection{Non-reversible dynamics}
To make a similar argument work in the non-reversible case seems a bit more hopeless, because there we cannot hope to show that the $\Gamma_n^\nu(\cdot)$ terms all vanish as this would imply reversibility, c.f. Proposition ~\ref{proposition:characterizations-reversibility}. 

In some sense, assuming reversibility allows to reduce it to a very local question, whereas mere stationarity is a global question. Compare this to Proposition 2.8 in \cite{liggett_interacting_2005} and the preceding discussion. On an intuitive level this can already be seen when considering continuous-time Markov chains on a finite state space. Under the assumption of reversibility, every edge of the transition graph is in equilibrium, whereas the weaker assumption of time-stationarity just implies that for every state the inflow and outflow of probability mass are equal. 

An alternative but related approach to show a result in the spirit of Corollary~\ref{corollary:no-periodic} for non-reversible interacting particle system would be to show that the relative entropy density as in e.g.~\cite{jahnel_dynamical_2022} is really a true Lyapunov function, in the sense that it is strictly negative for non-stationary $\nu$. 
However, working with non-shift-invariant measures requires to work with the $\limsup$ instead of the more convenient representations derived in \cite{jahnel_dynamical_2022} and additionally any argument of this type has to use the geometry of $d=1,2$ explicitly, because the above cannot be true in dimensions $d>2$ for short-range systems and in dimensions $d=1,2$ for long-range systems, as the examples in \cite{jahnel_class_2014} and \cite{jahnel_time_periodic_2024} show.

\section*{Acknowledgements}
The authors would like to thank the two anonymous referees for their insightful feedback that helped to find and fix a mistake in a previous version of this article. Additionally, the authors thank Aernout van Enter, Georg Menz, and Marek Biskup for comments and discussions that helped to clarify some technical details. 
Moreover, the authors acknowledge the financial support of the Leibniz Association within the Leibniz Junior Research Group on \textit{Probabilistic Methods for Dynamic Communication Networks} as part of the Leibniz Competition.

\bibliography{references}
\bibliographystyle{alpha}

\end{document}